\documentclass[reqno]{amsart}
\usepackage{amssymb,verbatim,mathtools}
\usepackage{tikz}
\usetikzlibrary{decorations.markings}
\usepackage[colorlinks=true,citecolor=blue,linkcolor=green,urlcolor=red]{hyperref}

\numberwithin{equation}{section}

\allowdisplaybreaks[4]

\theoremstyle{plain}
\newtheorem{thm}{Theorem}[section]
\newtheorem{lem}[thm]{Lemma}
\newtheorem{cor}[thm]{Corollary}
\newtheorem{problem}[thm]{Problem}
\theoremstyle{remark}
\newtheorem{rem}[thm]{Remark}

\DeclareMathOperator{\td}{d\!}
\DeclareMathOperator{\ti}{i}
\DeclareMathOperator{\te}{e}
\DeclareMathOperator{\bell}{B}
\DeclareMathOperator{\arctanh}{arctanh}
\DeclareMathOperator{\Res}{Res}
\DeclareMathOperator{\sech}{sech}
\newcommand{\fall}[2]{\langle #1\rangle_{#2}}

\begin{document}

\title[Wang's recursive sequence and monotonicity]
{Monotonicity of Wang's recursive sequence and logarithmically complete monotonicity of the reciprocal of the inverse tangent}

\author[F. Qi]{Feng Qi}
\address{School of Mathematics and Physics, Hulunbuir University, Hulunbuir 021008, Inner Mongolia, China;
17709 Sabal Court, University Place, Dallas, TX 75252-8024, USA}
\email{\href{mailto:qifeng618@gmail.com}{qifeng618@gmail.com}}
\urladdr{\url{https://orcid.org/0000-0001-6239-2968}}

\begin{abstract}
In this paper, we study a sequence \(b_n\), posed by Fei Wang, which is defined by \(b_1=\frac13\) and a convolution-type recursive relation. We derive the generating functions of \(nb_n\) and of \(b_n\), the latter of which is \(\ln\frac{\arctanh x}{x}\), and represent \(b_n\) as a moment integral over \([0,1]\) with a positive density. From this representation we prove that the sequence \(b_n\) is positive and logarithmically convex, that the sequence \(nb_n\) is strictly decreasing, and that the ratio \(\frac{b_{n+1}}{b_n}\) is strictly increasing in \(n\in\mathbb{N}\) with limit \(1\), which answers Wang's problem. We also give a determinantal expression and closed-form formulas for \(b_n\) involving Stirling numbers of the first kind and partial Bell polynomials. Finally, we prove that the function \(\frac1{\arctan x}\) is logarithmically completely monotonic on \((0,\infty)\).
\end{abstract}

\subjclass{Primary 11B83; Secondary 05A19, 26A48, 30E20, 33B10, 44A10}

\keywords{recursive relation; generating function; inverse hyperbolic tangent function; ratio of two consecutive terms; logarithmic convexity; integral representation; logarithmically completely monotonic function; partial Bell polynomial; Stirling number of the first kind; determinant}

\thanks{This paper was typeset using \AmS-\LaTeX}

\maketitle

\section{Introduction}

On 18 July 2025, Professor Fei Wang (Zhejiang Polytechnic University of Mechanical and Electrical Engineering, China) asked the author the following problem.

\begin{problem}\label{Wang-Fei-Problem}
For $n\in\mathbb{N}$, let $b_1=\frac{1}{3}$ and
\begin{equation}\label{Wang-Fei-Recur}
b_{n+1}=\frac{1}{n+1}\Biggl[\frac{n+1}{2n+3}-\sum_{k=0}^{n-1}\frac{(k+1)b_{k+1}}{2(n-k)+1}\Biggr].
\end{equation}
Prove that the sequence $T_n=\frac{b_{n+1}}{b_n}$ is increasing in $n\in\mathbb{N}$.
\end{problem}

This problem was posted on the website \url{https://mathoverflow.net/q/498467} on 30 July 2025 by the author.
\par
On 1 August 2025, Wang further told the author that the generating function of the sequence $b_n$ for $n\in\mathbb{N}$ is related to $\ln\frac{\arctan x}{x}$.
The precise statements are given in Theorem~\ref{arctanh-gen-bn-thm} and Remark~\ref{arctan-gen-bn-rem} below:
one has
\[\ln\frac{\arctanh x}{x}=\sum_{n=1}^{\infty}b_nx^{2n}
\quad\text{and}\quad
\ln\frac{\arctan x}{x}=\sum_{n=1}^{\infty}(-1)^nb_nx^{2n}.
\]
\par
Based on the above information, in this paper we study the sequence $b_n$ and its generating function.
The paper is organized as follows.
\begin{enumerate}
\item
In Section~\ref{sec-notation}, we fix notation and recall some basic facts.
\item
In Section~\ref{sec-gen-fun}, starting from the recursive relation~\eqref{Wang-Fei-Recur}, we derive the generating functions of the sequences $nb_n$ and $b_n$.
\item
In Section~\ref{sec-solution}, we represent $b_n$ as a moment of a positive density (Theorem~\ref{b(n)-int-rep-thm}),
and use this to solve Problem~\ref{Wang-Fei-Problem} in a stronger form (Theorem~\ref{Wang-Fei-solution-thm}).
\item
In Section~\ref{sec-closed}, we give closed-form expressions for $b_n$: a determinantal expression (Subsection~\ref{sec-determinant}), and formulas in terms of the Stirling numbers of the first kind and the partial Bell polynomials (Subsection~\ref{sec-Bell}).
\item
In Section~\ref{sec-logCM}, we prove that $\frac{1}{\arctan x}$ is logarithmically completely monotonic on $(0,\infty)$ and that $\frac{\arctanh x}{x}$ is logarithmically absolutely monotonic on $(0,1)$.
\end{enumerate}

\section{Notation and preliminaries}\label{sec-notation}

This section collects the notation and the elementary facts that are used throughout the paper: the standard sets of integers, the rising and falling factorials, the binomial coefficients, and the power series of $\arctanh x$ together with its relation to $\arctan x$.
\par
Throughout this paper, we use the notation $\mathbb{N}=\{1,2,\dotsc\}$ and $\mathbb{N}_0=\mathbb{N}\cup\{0\}$.
\par
For $\lambda\in\mathbb{C}$ and $\ell\in\mathbb{N}_0$, the rising factorial (the Pochhammer symbol, see~\cite[p.~256]{abram}) and the falling factorial are defined by
\begin{equation*}
(\lambda)_\ell=\prod_{k=0}^{\ell-1}(\lambda+k)
\quad\text{and}\quad
\fall{\lambda}{\ell}=\prod_{k=0}^{\ell-1}(\lambda-k),
\end{equation*}
with $(\lambda)_0=\fall{\lambda}{0}=1$, and $\binom{\lambda}{\ell}=\frac{\fall{\lambda}{\ell}}{\ell!}$ is the extended binomial coefficient (see~\cite[p.~154]{GKP-Concrete-Math-2nd} and~\cite[p.~17, Identity~17]{Spivey-art-2019}); in particular $\binom{m}{\ell}=0$ for integers $0\le m<\ell$.
\par
The symbol $\ln$ denotes the principal branch of the logarithm, and
\[
\arctanh z=\frac12\ln\frac{1+z}{1-z}, \quad z\in\mathbb{C}\setminus((-\infty,-1]\cup[1,\infty)).
\]
Recall the power series
\begin{equation*}
\arctanh x=\sum_{n=0}^{\infty}\frac{x^{2n+1}}{2n+1},\quad |x|<1,
\end{equation*}
see~\cite[p.~61]{Gradshteyn-Ryzhik-Table-8th}, and the relation
\begin{equation}\label{arctan=arctanh}
\arctan z=\frac{\arctanh (\ti z)}{\ti},
\end{equation}
see~\cite[p.~xxix]{Gradshteyn-Ryzhik-Table-8th}.

\section{Generating functions of the sequences \texorpdfstring{$nb_n$}{nb\_n} and \texorpdfstring{$b_n$}{b\_n}}\label{sec-gen-fun}

In this section, starting from the recursive relation~\eqref{Wang-Fei-Recur}, we derive the generating functions of the sequences $nb_n$ and $b_n$ for $n\in\mathbb{N}$ step by step.
\par
We first collect the analytic properties of the function
\begin{equation}\label{f(t)-dfn}
f(t)=\sum_{n=0}^{\infty}\frac{t^n}{2n+1}=\frac{\arctanh \sqrt{t}\,}{\sqrt{t}\,},\quad |t|<1,
\end{equation}
which will be used in both this section and Section~\ref{sec-solution}.

\begin{lem}\label{f-analytic-lem}
The function $f$ defined by~\eqref{f(t)-dfn} extends analytically to $\Omega=\mathbb{C}\setminus[1,\infty)$ and $\Re f(t)>0$ for all $t\in\Omega$.
In particular, $f$ has no zeros in $\Omega$, and $F=\ln f$ is analytic in $\Omega$, satisfies 
\[F(0)=0,\quad F(\bar t)=\overline{F(t)},\quad |\Im F(t)|<\frac{\pi}2.
\]
\end{lem}

\begin{proof}
Let $D=\mathbb{C}\setminus((-\infty,-1]\cup[1,\infty))$.
The function 
\[
\varphi(x)=\frac{\arctanh x}{x}, 
\]
with $\varphi(0)=1$, is analytic and even on $D$, and $x\in D$ if and only if $x^2\in\Omega$.
Hence $f(t)=\varphi\bigl(\sqrt{t}\,\bigr)$ is a well defined analytic function on $\Omega$ that agrees with~\eqref{f(t)-dfn} in $|t|<1$.
\par
If $t\in(-\infty,1)$, then $f(t)>0$: for $t\in[0,1)$ by the series~\eqref{f(t)-dfn}, and for $t<0$ because $f(t)=\frac{\arctan\sqrt{-t}\,}{\sqrt{-t}\,}$ by~\eqref{arctan=arctanh}.
\par
Let $\Im t>0$ and let $x=\sqrt{t}\,$ be the principal root, so that $0<\arg x<\frac\pi2$, and put $z=\frac{1+x}{1-x}$.
Then
\[
|1+x|^2-|1-x|^2=4\Re x>0\quad \text{and}\quad \Im z=\frac{2\Im x}{|1-x|^2}>0,
\] 
so that $|z|>1$ and $0<\arg z<\pi$.
Hence $\arctanh x=\frac12\ln z$ lies in the open first quadrant, that is, $0<\arg\arctanh x<\frac\pi2$.
Consequently 
\[
\arg f(t)=\arg\arctanh x-\arg x\in\biggl(-\frac\pi2,\frac\pi2\biggr),
\] 
that is, $\Re f(t)>0$.
\par
The case $\Im t<0$ follows from $f(\bar t)=\overline{f(t)}$.
\par
The remaining assertions are immediate consequences of $\Re f>0$.
\end{proof}

\begin{thm}
The generating function of the sequence $nb_n$ for $n\in\mathbb{N}$ is
\begin{equation}\label{nbn-gen-funct-eq}
\frac{1}{2}\biggl(\frac{1}{1-x^2}\frac{x}{\arctanh x}-1\biggr)
=\sum_{n=1}^{\infty}nb_nx^{2n}, \quad |x|<1,
\end{equation}
where $\frac{x}{\arctanh x}=1$ at $x=0$.
\end{thm}

\begin{proof}
Let $a_n=nb_n$ for $n\in\mathbb{N}$, so that $a_1=\frac13$. Multiplying the recursive relation~\eqref{Wang-Fei-Recur} by $n+1$ and shifting the summation index, we obtain
\begin{equation*}
a_{n+1}=\frac{n+1}{2n+3}-\sum_{k=0}^{n-1}\frac{a_{k+1}}{2(n-k)+1}
=\frac{n+1}{2n+3}-\sum_{k=1}^{n}\frac{a_{k}}{2(n-k)+3}, \quad n\in\mathbb{N}.
\end{equation*}
Multiplying by $t^{n+1}$ and summing over $n\in\mathbb{N}$, we obtain, as an identity between formal power series,
\begin{equation*}
\sum_{n=1}^{\infty}a_{n+1}t^{n+1}=\sum_{n=1}^{\infty}\frac{n+1}{2n+3}t^{n+1} -\sum_{n=1}^{\infty}\Biggl[\sum_{k=1}^{n}\frac{a_{k}}{2(n-k)+3}\Biggr]t^{n+1},
\end{equation*}
where
\begin{gather*}
\sum_{n=1}^{\infty}a_{n+1}t^{n+1}=\sum_{n=1}^{\infty}a_{n}t^n-a_1t
=A(t)-\frac{1}{3}t,\\
\sum_{n=1}^{\infty}\frac{n+1}{2n+3}t^{n+1}=\sum_{n=2}^{\infty}\frac{n}{2n+1}t^{n}
=B(t)-\frac{1}{3}t,
\intertext{and}
\sum_{n=1}^{\infty}\Biggl[\sum_{k=1}^{n}\frac{a_{k}}{2(n-k)+3}\Biggr]t^{n+1}
=t\Biggl(\sum_{n=1}^{\infty}a_{n}t^n\Biggr)\Biggl(\sum_{n=0}^{\infty}\frac{t^n}{2n+3}\Biggr)
=tA(t)C(t),
\end{gather*}
with
\begin{equation*}
A(t)=\sum_{n=1}^{\infty}a_{n}t^n, \quad B(t)=\sum_{n=1}^{\infty}\frac{n}{2n+1}t^{n}, \quad C(t)=\sum_{n=0}^{\infty}\frac{t^n}{2n+3}.
\end{equation*}
Therefore, 
\[
A(t)-\frac{1}{3}t=B(t)-\frac{1}{3}t-tA(t)C(t), 
\]
that is,
\begin{equation}\label{gen-func-eq-A(x)}
A(t)\bigl[1+tC(t)\bigr]=B(t).
\end{equation}
\par
Since 
\[
\frac{n}{2n+1}=\frac12\biggl(1-\frac1{2n+1}\biggr)\quad \text{and}\quad tC(t)=\sum_{n=0}^{\infty}\frac{t^{n+1}}{2n+3}=f(t)-1,
\] 
where $f$ is defined by~\eqref{f(t)-dfn}, we have
\begin{equation*}
1+tC(t)=f(t)
\quad\text{and}\quad
B(t)=\frac12\Biggl(\sum_{n=0}^{\infty}t^n-\sum_{n=0}^{\infty}\frac{t^n}{2n+1}\Biggr)=\frac12\biggl[\frac{1}{1-t}-f(t)\biggr]
\end{equation*}
for $|t|<1$.
By Lemma~\ref{f-analytic-lem}, the function $f$ has no zero in the unit disc, so $\frac{B}{f}$ is analytic in $|t|<1$ and its Maclaurin series is the formal quotient given by~\eqref{gen-func-eq-A(x)}.
Hence
\begin{equation*}
A(t)=\frac{B(t)}{f(t)}=\frac{1}{2}\biggl[\frac{1}{1-t}\frac{1}{f(t)}-1\biggr]=\sum_{n=1}^{\infty}nb_nt^n, \quad |t|<1.
\end{equation*}
Replacing $t$ by $x^2$ and using $f(x^2)=\frac{\arctanh x}{x}$ gives~\eqref{nbn-gen-funct-eq}.
\end{proof}

\begin{rem}
Put $d_n=nb_n-(n-1)b_{n-1}$ for $n\ge2$.
The generating function~\eqref{nbn-gen-funct-eq} can be reformulated as
\begin{equation*}
\frac{x}{\arctanh x}
=1-\frac{1}{3}x^{2}+2\sum_{n=2}^{\infty}d_nx^{2n}, \quad |x|<1.
\end{equation*}
Replacing $x$ by $\ti x$ and using the relation~\eqref{arctan=arctanh}, we derive
\begin{equation}\label{arctan-recipr}
\frac{x}{\arctan x}
=1+\frac{1}{3}x^{2}+2\sum_{n=2}^{\infty}(-1)^nd_nx^{2n}, \quad |x|<1.
\end{equation}
In the same way, dividing~\eqref{nbn-gen-funct-eq} by $\frac x2$ and replacing $x$ by $\ti x$ yield
\begin{equation*}
\frac{1}{(1+x^2)\arctan x}=\frac1x+2\sum_{n=1}^{\infty}(-1)^nnb_nx^{2n-1}, \quad 0<|x|<1.
\end{equation*}
\end{rem}

\begin{thm}\label{arctanh-gen-bn-thm}
The generating function of the sequence $b_n$ for $n\in\mathbb{N}$ is
\begin{equation}\label{b-n-gen-f}
\ln\frac{\arctanh x}{x}=\sum_{n=1}^{\infty}b_nx^{2n}, \quad |x|<1.
\end{equation}
\end{thm}

\begin{proof}
Dividing~\eqref{nbn-gen-funct-eq} by $\frac x2$, we can reformulate it as
\begin{equation}\label{nbn-gen-funct-reform}
\frac{1}{1-x^2}\frac{1}{\arctanh x}-\frac1x
=\sum_{n=1}^{\infty}2nb_nx^{2n-1}, \quad 0<|x|<1.
\end{equation}
Integrating both sides of~\eqref{nbn-gen-funct-reform} over the segment from $0$ to $s\in(-1,1)$ and using
\[
\bigl(\ln\arctanh x\bigr)'=\frac{1}{(1-x^2)\arctanh x}\quad \text{and}\quad \lim_{x\to0}\frac{\arctanh x}{x}=1
\] 
give
\begin{equation*}
\int_{0}^{s}\biggl(\frac{1}{1-x^2}\frac{1}{\arctanh x}-\frac1x\biggr)\td x
=\ln\frac{\arctanh s}{s}
=\sum_{n=1}^{\infty}b_ns^{2n}, \quad |s|<1.
\end{equation*}
This is~\eqref{b-n-gen-f}.
\end{proof}

\begin{rem}\label{arctan-gen-bn-rem}
Replacing $x$ by $\ti x$ in~\eqref{b-n-gen-f} and using the relation~\eqref{arctan=arctanh}, we obtain
\begin{equation}\label{b-n-gen-f-arctan}
\ln\frac{\arctan x}{x}=\sum_{n=1}^{\infty}(-1)^nb_nx^{2n}, \quad |x|<1.
\end{equation}
In terms of the function $f$ in~\eqref{f(t)-dfn}, Theorem~\ref{arctanh-gen-bn-thm} says that
\[
F(t)=\ln f(t)=\sum_{n=1}^{\infty}b_nt^n
\]
for $|t|<1$.
\end{rem}

\section{Solution of Problem~\ref{Wang-Fei-Problem}}\label{sec-solution}

In this section, we solve Problem~\ref{Wang-Fei-Problem} by means of the generating function obtained in Section~\ref{sec-gen-fun}. We first represent $b_n$ as an integral (Theorem~\ref{b(n)-int-rep-thm}), then deduce the monotonicity statements asked for in the problem (Theorem~\ref{Wang-Fei-solution-thm}), and finally record the consequences for the sequence $nb_n$ (Corollary~\ref{b(n)-positive-cor}).
\par
The idea is to show that $b_n$ is a moment sequence of a positive density. Such a sequence is automatically logarithmically convex by the Cauchy--Schwarz inequality, which is exactly the statement that $\frac{b_{n+1}}{b_n}$ is increasing.

\begin{thm}\label{b(n)-int-rep-thm}
Let
\begin{equation*}
w(s)=\frac1\pi\arctan\frac{\pi}{2\arctanh \sqrt{s}\,},\quad s\in(0,1).
\end{equation*}
Then $w$ is positive and strictly decreasing on $(0,1)$, with $w(0^+)=\frac12$ and $w(1^-)=0$, and for every $n\in\mathbb{N}$,
\begin{align}
b_n&=\int_0^1s^{n-1}w(s)\td s\label{b(n)-moment}\\
&=\frac{2}{n}\int_0^\infty\frac{\tanh^{2n}y}{\pi^2+4y^2}\td y.\label{b(n)-tanh}
\end{align}
\end{thm}

\begin{proof}
The monotonicity and limits of $w$ are clear, because $\arctanh \sqrt{s}\,$ increases from $0$ to $\infty$ on $(0,1)$.
\par
\paragraph{\indent\it Step~1: boundary values.}
Let $f$ and $F=\ln f$ be as in Lemma~\ref{f-analytic-lem}, so that 
\[
F(t)=\sum_{m=1}^{\infty}b_mt^m, \quad |t|<1
\] 
by Theorem~\ref{arctanh-gen-bn-thm}.
For $t>1$ and $\varepsilon\downarrow0$, the point $x=\sqrt{t+\ti\varepsilon}\,$ tends to $\sqrt{t}\,>1$, and $z=\frac{1+x}{1-x}$ tends to the negative number $-\frac{\sqrt{t}\,+1}{\sqrt{t}\,-1}$ from the upper half plane, so 
\[
\ln z\to\ln\frac{\sqrt{t}\,+1}{\sqrt{t}\,-1}+\ti\pi.
\]
Hence 
\[
\arctanh x\to L(t)+\frac{\ti\pi}{2}
\] 
with
\begin{equation*}
L(t)=\frac12\ln\frac{\sqrt{t}\,+1}{\sqrt{t}\,-1}=\arctanh \frac{1}{\sqrt{t}\,}>0,
\end{equation*}
and, since $\sqrt{t}\,>0$,
\begin{gather}\label{Im-F-boundary}
\Im F(t+\ti0)=\arg\Bigl(L(t)+\frac{\ti\pi}2\Bigr)=\arctan\frac{\pi}{2L(t)}\in\Bigl(0,\frac\pi2\Bigr),\\
F(t-\ti0)=\overline{F(t+\ti0)}.\notag
\end{gather}
\par
\paragraph{\indent\it Step~2: growth estimates.}
Let $|t|\ge25$, $\Im t\ge0$, and $x=\sqrt{t}\,$ (so $\Re x\ge0$, $\Im x\ge0$, and $|x|\ge5$).
Then 
\[
z=\frac{1+x}{1-x}=-1-\frac{2}{x-1}
\] 
lies in the disc of radius $\frac12$ centred at $-1$, so $\frac{5\pi}{6}\le\arg z\le\pi$ and, because $\Re x\ge0$ gives $|z|\ge1$, also
\[
0\le\ln|z|\le\ln\frac32<\frac12.
\]
Therefore
\[
1<\frac{5\pi}{12}\le\Im\arctanh x\le|\arctanh x|\le\frac{\ln|z|+\pi}2<2,
\] 
i.e. 
\[
\frac1{|x|}<|f(t)|<\frac2{|x|}.
\]
Hence 
\[
|\ln|f(t)||\le\frac12\ln|t|,\quad |t|\ge25,
\] 
and since $|F|\le|\ln|f||+\frac\pi2$ and $F(\bar t)=\overline{F(t)}$, we obtain
\[
|F(t)|\le\frac12\ln|t|+3,\quad |t|\ge25
\] 
in the whole plane.
\par
Let $t\in\Omega$ with $|t-1|$ small and $x=\sqrt{t}\,$ in the closed right half plane.
Then $|1-x|\le|t-1|$ and $|1+x|\ge1$, so 
\[
|\arctanh x|\ge\frac12\ln|z|\ge\frac12\ln\frac{1}{|t-1|},
\] 
while 
\[
|\arctanh x|\le\frac12\biggl(\ln\frac{9}{|t-1|}+\pi\biggr).
\]
Hence 
\[
|F(t)|\le\ln\ln\frac1{|t-1|}+C
\] 
near $t=1$ for some constant $C$.
\par
\paragraph{\indent\it Step~3: the contour integral.}
Fix $n\in\mathbb{N}$.
Since $F(t)t^{-n-1}$ has the residue $b_n$ at $t=0$ and is analytic in $\Omega\setminus\{0\}$, Cauchy's residue theorem applied to the following positively oriented closed contour gives $2\pi\ti b_n$ as the sum of four integrals:
the upper edge of the cut from $1+\delta$ to $R$; the circle $|t|=R$ counterclockwise from $R+\ti0$ to $R-\ti0$;
the lower edge of the cut from $R$ back to $1+\delta$; and the circle $|t-1|=\delta$ clockwise around $1$.
\par
By Step~2, the integral over $|t|=R$ is bounded by 
\[
2\pi R\biggl(\frac12\ln R+3\biggr)R^{-n-1}\to0,\quad R\to\infty,
\] 
because $n\ge1$.
\par
A standard estimate of the integral over the small circle $|t-1|=\delta$ is as follows:
\begin{enumerate}
\item
Length of the circle: at most $2\pi\delta$.
\item
Bound on $F$: by the second paragraph of Step~2, there is a constant $C$ such that
\[
|F(t)|\le\ln\ln\frac1{|t-1|}+C=\ln\ln\frac1\delta+C
\] 
on that circle.
\item
Bound on the kernel: for $\delta\le\frac12$, we have
\[
|t^{-n-1}|\le(1-\delta)^{-n-1}\le2^{n+1}.
\]
\end{enumerate}
Multiplying these gives a bound of 
\[
2\pi\cdot2^{n+1}\delta\biggl(\ln\ln\frac1\delta+C\biggr).
\] 
Since $\ln\ln\frac1\delta\to\infty$ as $\delta\to0^+$, the additive $C$ can be absorbed into the main term for small $\delta$. The whole bound is then at most $C'\delta\ln\ln\frac1\delta$, for example with $C'=2\pi\cdot2^{n+1}(1+C)$ once $\ln\ln\frac1\delta\ge1$. Because $\delta\ln\ln\frac1\delta\to0$, this integral vanishes in the limit.
\par
By~\eqref{Im-F-boundary}, the two edges of the cut contribute in the limit
\begin{equation*}
2\pi\ti b_n=\int_1^\infty[F(t+\ti0)-F(t-\ti0)]t^{-n-1}\td t=2\ti\int_1^\infty\Im F(t+\ti0)t^{-n-1}\td t,
\end{equation*}
and the last integral converges absolutely, because $0<\Im F(t+\ti0)<\frac\pi2$.
Thus
\begin{equation*}
b_n=\frac1\pi\int_1^\infty\arctan\frac{\pi}{2L(t)}t^{-n-1}\td t.
\end{equation*}
Substituting $t=\frac1s$ and noting $L\bigl(\frac1s\bigr)=\arctanh \sqrt{s}\,$ yield~\eqref{b(n)-moment}.
\par
\paragraph{\indent\it Step~4: the formula~\eqref{b(n)-tanh}.}
Substituting $s=\tanh^2y$ in~\eqref{b(n)-moment}, so that
\[
\arctanh \sqrt{s}=y\quad \text{and}\quad \td s=2\tanh y\sech^2y\td y,
\]
gives
\begin{equation*}
b_n=\frac2\pi\int_0^\infty\tanh^{2n-1}y\sech^2y\arctan\frac{\pi}{2y}\td y.
\end{equation*}
Integrating by parts with 
\[
\operatorname{d}\biggl(\frac{\tanh^{2n}y}{2n}\biggr)=\tanh^{2n-1}y\sech^2y\td y\quad \text{and} \quad \frac{\td}{\td y}\arctan\frac{\pi}{2y}=-\frac{2\pi}{\pi^2+4y^2},
\] 
where the boundary terms vanish, we arrive at~\eqref{b(n)-tanh}.
\end{proof}

\begin{thm}[Solution of Problem~\ref{Wang-Fei-Problem}]\label{Wang-Fei-solution-thm}
For every $n\in\mathbb{N}$, we have
\begin{equation}\label{log-convex-ineq}
b_nb_{n+2}>b_{n+1}^2.
\end{equation}
Consequently, the sequence $T_n=\frac{b_{n+1}}{b_n}$ is strictly increasing in $n\in\mathbb{N}$.
Moreover, 
\[
T_1=\frac{13}{30},\quad 0<T_n<1,\quad \text{and}
\quad \lim_{n\to\infty}T_n=1.
\]
\end{thm}

\begin{proof}
By~\eqref{b(n)-moment}, it is clear that 
\[
b_n=\int_0^1s^{n-1}w(s)\td s>0
\] 
for all $n\in\mathbb{N}$.
\par
Writing $s^n=s^{(n-1)/2}s^{(n+1)/2}$ and applying the Cauchy--Schwarz inequality, we obtain
\begin{multline*}
b_{n+1}^2=\biggl(\int_0^1\bigl[s^{(n-1)/2}\sqrt{w(s)}\,\bigr]\bigl[s^{(n+1)/2}\sqrt{w(s)}\,\bigr]\td s\biggr)^2
\\
\le\int_0^1s^{n-1}w(s)\td s\int_0^1s^{n+1}w(s)\td s=b_nb_{n+2}.
\end{multline*}
Equality would require the two functions $s^{(n-1)/2}\sqrt{w(s)}\,$ and $s^{(n+1)/2}\sqrt{w(s)}\,$ to be proportional, which is impossible since their ratio $s$ is not constant and $w>0$. This proves~\eqref{log-convex-ineq}.
\par
Dividing~\eqref{log-convex-ineq} by $b_nb_{n+1}>0$ gives $T_n<T_{n+1}$.
\par
Since $0<s<1$, we have 
\[
b_{n+1}=\int_0^1s^nw(s)\td s<b_n,
\] 
that is, $0<T_n<1$. 
\par
From the recursive relation~\eqref{Wang-Fei-Recur}, it is obvious that 
\[
b_2=\frac{1}{2}\biggl(\frac25-\frac19\biggr)=\frac{13}{90},
\] 
hence $T_1=\frac{13}{30}$.
\par
Finally, let $\varepsilon\in(0,1)$. For $n\in\mathbb{N}$,
\begin{multline*}
b_{n+1}\ge\int_{1-\varepsilon}^1s^nw(s)\td s\ge(1-\varepsilon)\int_{1-\varepsilon}^1s^{n-1}w(s)\td s
\\
=(1-\varepsilon)\biggl[b_n-\int_0^{1-\varepsilon}s^{n-1}w(s)\td s\biggr].
\end{multline*}
On the other hand, with $c_\varepsilon=\int_{1-\varepsilon/2}^1w(s)\td s>0$, we have
\[
b_n\ge\int_{1-\varepsilon/2}^1s^{n-1}w(s)\td s\ge c_\varepsilon\Bigl(1-\frac\varepsilon2\Bigr)^{n-1}
\]
and
\[
\int_0^{1-\varepsilon}s^{n-1}w(s)\td s\le\frac12(1-\varepsilon)^{n-1}.
\]
Hence 
\[
\frac{1}{b_n}\int_0^{1-\varepsilon}s^{n-1}w(s)\td s\to0, \quad n\to\infty
\] 
and $\liminf_{n\to\infty}T_n\ge1-\varepsilon$.
Letting $\varepsilon\to0^+$ and using $0<T_n<1$ give the limit $\lim_{n\to\infty}T_n=1$.
\end{proof}

\begin{cor}\label{b(n)-positive-cor}
For $n\in\mathbb{N}$, we have $b_n>0$. The sequence $nb_n$ is positive, strictly decreasing, and logarithmically convex, because
\begin{equation}\label{nbn-tanh}
nb_n=2\int_0^\infty\frac{\tanh^{2n}y}{\pi^2+4y^2}\td y.
\end{equation}
In particular, $d_n=nb_n-(n-1)b_{n-1}<0$ for $n\ge2$.
\end{cor}

\begin{proof}
The formula~\eqref{nbn-tanh} is precisely~\eqref{b(n)-tanh}.
\par
Since $0<\tanh y<1$ for $y>0$, we have $\tanh^{2n+2}y<\tanh^{2n}y$, which gives the strict decrease.
\par
Writing 
\[
\tanh^{2n+2}y=\tanh^{n}y\tanh^{n+2}y
\] 
and applying the Cauchy--Schwarz inequality with respect to the positive measure $\frac{\td y}{\pi^2+4y^2}$, we obtain 
\[
[(n+1)b_{n+1}]^2<(nb_n)[(n+2)b_{n+2}], 
\]
the strictness following because $\tanh^ny$ and $\tanh^{n+2}y$ are not proportional.
\end{proof}

\begin{rem}
Numerically, we have
\[
T_1=0.433\dotsc,\quad T_2=\frac{251\cdot90}{2835\cdot13}=0.612\dotsc,\quad T_3=0.707\dotsc,\quad T_4=0.765\dotsc,
\] 
in agreement with Theorem~\ref{Wang-Fei-solution-thm}.
\par
The first values 
\[
b_1=\frac13,\quad b_2=\frac{13}{90},\quad b_3=\frac{251}{2835},\quad b_4=\frac{3551}{56700},\quad b_5=\frac{22417}{467775}
\] 
follow from~\eqref{Wang-Fei-Recur} and from~\eqref{b(n)-moment}.
\end{rem}

\section{Closed-form expressions of the sequence \texorpdfstring{$b_n$}{b\_n}}\label{sec-closed}

This section gives two kinds of closed-form expressions for $b_n$, which complement the integral representation of Section~\ref{sec-solution}. In Subsection~\ref{sec-determinant}, we write $b_n$ as a determinant of order $n$ with explicit entries. In Subsection~\ref{sec-Bell}, we express $b_n$ by finite sums involving the Stirling numbers of the first kind and the partial Bell polynomials.

\subsection{Determinantal expression}\label{sec-determinant}

In this subsection, we express $b_n$ as a determinant by applying a formula of Bourbaki for the derivatives of a quotient of two functions to the logarithmic derivative $\frac{f'(t)}{f(t)}$ of the function $f(t)$ in~\eqref{f(t)-dfn}. We also give a second proof of the positivity of $b_n$ based on a criterion for the signs of power series coefficients.

\begin{lem}[{\cite[p.~40, Exercise~5]{Bourbaki-Spain-2004}}]\label{Bourbaki-lem}
For $n\in\mathbb{N}_0$, let $u(x)$ and $v(x)\ne0$ be $n$ times differentiable functions, and let
\begin{equation*}
W_{(n+1)\times(n+1)}(x)=\begin{pmatrix}U_{(n+1)\times1}(x) & V_{(n+1)\times n}(x)\end{pmatrix}_{(n+1)\times(n+1)}
\end{equation*}
be the $(n+1)\times(n+1)$ matrix whose first column has the elements $u_{k,1}(x)=u^{(k-1)}(x)$ for $1\le k\le n+1$ and whose other columns have the elements $v_{\ell,j}(x)=\binom{\ell-1}{j-1}v^{(\ell-j)}(x)$ for $1\le \ell\le n+1$ and $1\le j\le n$.
Then
\begin{equation*}
\frac{\operatorname{d}^n}{\td x^n}\biggl[\frac{u(x)}{v(x)}\biggr]
=(-1)^n\frac{\bigl|W_{(n+1)\times(n+1)}(x)\bigr|}{v^{n+1}(x)}, \quad n\in\mathbb{N}_0.
\end{equation*}
\end{lem}

\begin{thm}
For $n\in\mathbb{N}$, the sequence $b_n$ has the determinantal expression
\begin{equation}\label{b(n)-determ}
\begin{aligned}
b_{n}&=\frac{(-1)^{n+1}}{n!}
\begin{vmatrix}
\frac{1!}{3} & \frac{0!}{1} & 0 & \dotsm & 0 & 0& 0\\
\frac{2!}{5} & \frac{1!}{3} & \frac{0!}{1} & \dotsm&0 & 0& 0\\
\frac{3!}{7} & \frac{2!}{5} & \binom{2}1\frac{1!}{3} & \dotsm&0 & 0& 0\\
\vdots & \vdots & \vdots & \ddots & \vdots & \vdots & \vdots\\
\frac{(n-2)!}{2n-3} & \frac{(n-3)!}{2n-5} & \binom{n-3}1\frac{(n-4)!}{2n-7} & \dotsm & \binom{n-3}{n-4}\frac{1!}{3} & \frac{0!}{1} & 0\\
\frac{(n-1)!}{2n-1} & \frac{(n-2)!}{2n-3} & \binom{n-2}1\frac{(n-3)!}{2n-5} & \dotsm& \binom{n-2}{n-4}\frac{2!}{5} & \binom{n-2}{n-3}\frac{1!}{3} & \frac{0!}{1}\\
\frac{n!}{2n+1} & \frac{(n-1)!}{2n-1} & \binom{n-1}1\frac{(n-2)!}{2n-3} & \dotsm& \binom{n-1}{n-4}\frac{3!}{7} & \binom{n-1}{n-3}\frac{2!}{5} & \binom{n-1}{n-2}\frac{1!}{3}
\end{vmatrix}\\
&=\frac{(-1)^{n+1}}{n!}|e_{i,j}|_{n\times n},
\end{aligned}
\end{equation}
where
\begin{equation*}
e_{i,j}=\begin{dcases}
\frac{i!}{2i+1}, & j=1;\\
\binom{i-1}{j-2}\frac{(i-j+1)!}{2(i-j)+3}, & 2\le j\le i+1;\\
0, & j>i+1
\end{dcases}
\end{equation*}
for $i,j\in\mathbb{N}$ with $i,j\le n$.
\end{thm}

\begin{proof}
By Theorem~\ref{arctanh-gen-bn-thm} with $x^2=t$, we have 
\[
F(t)=\ln f(t)=\sum_{n=1}^{\infty}b_nt^n, \quad |t|<1,
\] 
where $f(t)=\sum_{n=0}^{\infty}\frac{t^n}{2n+1}$ is the function in~\eqref{f(t)-dfn}. Differentiating yields
\begin{equation}\label{(n+1)b(n+1)}
\frac{f'(t)}{f(t)}=\frac{\sum_{k=0}^{\infty}\frac{k+1}{2k+3}t^{k}}{\sum_{k=0}^{\infty}\frac{1}{2k+1}t^k}
=\sum_{k=0}^{\infty}(k+1)b_{k+1}t^{k}, \quad |t|<1.
\end{equation}
Therefore, 
\[
(m+1)!b_{m+1}=\frac{\operatorname{d}^m}{\td t^m}\biggl[\frac{f'(t)}{f(t)}\biggr]_{t=0},\quad m\in\mathbb{N}_0.
\]
\par
We apply Lemma~\ref{Bourbaki-lem} with $u=f'$ and $v=f$, and evaluate at $t=0$.
Since $f(0)=1$ and
\begin{equation*}
f^{(k)}(0)=\frac{k!}{2k+1},\quad f^{(k+1)}(0)=\frac{(k+1)!}{2k+3},\quad k\in\mathbb{N}_0,
\end{equation*}
the first column of $W_{(m+1)\times(m+1)}(0)$ has the elements $u^{(i-1)}(0)=\frac{i!}{2i+1}=e_{i,1}$ for $1\le i\le m+1$, and the element in row $i$ and column $j+1$ is
\begin{equation*}
\binom{i-1}{j-1}v^{(i-j)}(0)=\binom{i-1}{j-1}\frac{(i-j)!}{2(i-j)+1}=e_{i,j+1},\quad 1\le j\le m,
\end{equation*}
which vanishes for $j>i$ because $\binom{i-1}{j-1}=0$ there. Hence
\begin{equation*}
(m+1)!b_{m+1}=(-1)^m|e_{i,j}|_{(m+1)\times(m+1)},\quad m\in\mathbb{N}_0.
\end{equation*}
Replacing $m+1$ by $n$ yields~\eqref{b(n)-determ}.
\end{proof}

We recall a criterion on signs of power series coefficients, which gives another proof of the positivity of $b_n$ obtained in Corollary~\ref{b(n)-positive-cor}.

\begin{lem}[{\cite[Theorems~1 and~2]{Jurkat-PAMS-1954}}]\label{Theorem1MR64890}
Let
\begin{equation*}
k(t)=\sum_{\ell=0}^{\infty}k_\ell t^\ell,\quad q(t)=\sum_{\ell=0}^{\infty}q_\ell t^\ell,\quad p(t)=\sum_{\ell=0}^{\infty}p_\ell t^\ell
\end{equation*}
be formal series such that $k(t)=\frac{q(t)}{p(t)}$ and $p_\ell>0$ for $\ell\in\mathbb{N}_0$.
\begin{enumerate}
\item
If both of the sequences $\frac{p_{\ell+1}}{p_\ell}$ and $\frac{q_\ell}{p_\ell}$ are increasing in $\ell\in\mathbb{N}_0$, then $k_\ell\ge0$ for $\ell\in\mathbb{N}$.
\item
If the sequence $\frac{p_{\ell+1}}{p_\ell}$ is increasing in $\ell\in\mathbb{N}_0$ and the sequence $\frac{q_\ell}{p_\ell}$ is decreasing in $\ell\in\mathbb{N}_0$, then $k_\ell\le0$ for $\ell\in\mathbb{N}$.
\end{enumerate}
\end{lem}

\begin{rem}
Apply part~(1) of Lemma~\ref{Theorem1MR64890} to~\eqref{(n+1)b(n+1)} with $p_\ell=\frac{1}{2\ell+1}$ and $q_\ell=\frac{\ell+1}{2\ell+3}$.
The sequences 
\[
\frac{p_{\ell+1}}{p_\ell}=\frac{2\ell+1}{2\ell+3}\quad \text{and}\quad \frac{q_\ell}{p_\ell}=\frac{(\ell+1)(2\ell+1)}{2\ell+3}
\] 
are increasing in $\ell\in\mathbb{N}_0$, so $(\ell+1)b_{\ell+1}\ge0$ for $\ell\in\mathbb{N}$.
Together with $b_1=\frac13>0$, this shows again that $b_n\ge0$ for $n\in\mathbb{N}$. Corollary~\ref{b(n)-positive-cor} gives the strict inequality.
\end{rem}

\subsection{Explicit formulas via Stirling numbers and Bell polynomials}\label{sec-Bell}

In this subsection, we derive explicit formulas for $b_n$ from the Maclaurin series of powers of $\frac{\arctan z}{z}$ due to Grishaev and Sazonov. We first recall the necessary definitions and known results on the Bell polynomials and the Stirling numbers.
\par
In~\cite[p.~134, Theorem~A]{Comtet-Combinatorics-74}, the Bell polynomials of the second kind, also known as the partial Bell polynomials, are defined by
\begin{equation*}
\frac{\bell_{n,k}(x_1,x_2,\dotsc,x_{n-k+1})}{n!}
=\sum_{\substack{\ell_1,\ell_2,\dotsc,\ell_{n-k+1}\in\mathbb{N}_0\\ \sum_{i=1}^{n-k+1}i\ell_i=n\\
\sum_{i=1}^{n-k+1}\ell_i=k}} \prod_{i=1}^{n-k+1}\biggl[\frac{1}{\ell_i!}\biggl(\frac{x_i}{i!}\biggr)^{\ell_i}\biggr]
\end{equation*}
for $n,k\in\mathbb{N}_0$ such that $n\ge k$, with the special cases $\bell_{0,0}=1$ and 
\[
\bell_{n,0}(x_1,x_2,\dotsc,x_{n+1})=0
\] 
for $n\in\mathbb{N}$; by convention $\bell_{n,k}=0$ for $k>n$.
They have the generating function
\begin{equation}\label{113-final-formula}
\frac1{k!}\Biggl(\sum_{m=1}^{\infty} x_m\frac{t^m}{m!}\Biggr)^k
=\sum_{n=k}^{\infty} \bell_{n,k}(x_1,x_2,\dotsc,x_{n-k+1})\frac{t^n}{n!},\quad k\in\mathbb{N}_0,
\end{equation}
see~\cite[p.~133]{Comtet-Combinatorics-74}, and satisfy the Fa\`a di Bruno formula 
\begin{equation}\label{Bruno-Bell-Polynomial}
[f(g(z))]^{(k)}=\sum_{j=0}^k f^{(j)}(g(z)) \bell_{k,j}\bigl(g'(z),g''(z),\dotsc,g^{(k-j+1)}(z)\bigr),\quad k\in\mathbb{N}_0,
\end{equation}
where $f$ and $g$ are $k$ times differentiable functions, see~\cite{Ion-Chitescu-B-2017} and~\cite[p.~139, Theorem~C]{Comtet-Combinatorics-74}.
\par
The Stirling numbers of the first kind $s(n,k)$ for $n,k\in\mathbb{N}_0$ are generated~\cite[Theorem~3.14]{Mansour-Schork-B2016} by
\begin{equation*}
\biggl[\frac{\ln(1+z)}{z}\biggr]^k=\sum_{n=0}^\infty \frac{s(n+k,k)}{\binom{n+k}{k}}\frac{z^{n}}{n!}, \quad |z|<1.
\end{equation*}
In~\cite{Grishaev-Sazonov-2024-ITSF}, Grishaev and Sazonov established the Maclaurin series expansions
\begin{equation}\label{Grishaev-Sazonov-ITSF-1}
\biggl(\frac{\arctan z}{z}\biggr)^m=m!\sum_{k=0}^{\infty}(-1)^k G(m,k)z^{2k}
\end{equation}
and
\begin{equation}\label{Grishaev-Sazonov-ITSF-2}
\biggl(\frac{\arctan z}{z}\biggr)^r =1+\sum_{m=1}^{\infty}(-1)^m\Biggl[\sum_{k=1}^{2m}\frac{(-r)_k}{k!} \sum_{j=1}^{k}(-1)^j\binom{k}{j}j!G(j,m)\Biggr]z^{2m}
\end{equation}
for $m\in\mathbb{N}$, $r\in\mathbb{R}$, and $|z|<1$, where
\begin{equation*}
G(m,k)=\sum_{\ell=0}^{2k}\binom{2k+m-1}{2k-\ell}2^\ell\frac{s(\ell+m,m)}{(\ell+m)!}
\end{equation*}
for $m,k\in\mathbb{N}_0$. For $m=0$, we have $G(0,k)=\binom{2k-1}{2k}$, which equals $1$ for $k=0$ and $0$ for $k\in\mathbb{N}$, since $s(\ell,0)=0$ for $\ell\in\mathbb{N}$; thus~\eqref{Grishaev-Sazonov-ITSF-1} also holds for $m=0$.
See also~\cite[Remark~12]{series-arccos-up-42.tex} and~\cite[Section~1.3]{arcsine-liu-qi-mdpi.tex}.
\par
In this subsection, we establish a combinatorial identity for $G(m,k)$, a closed-form formula for specific values of the Bell polynomials $\bell_{2n,k}$, a closed-form formula for $b_n$, and a formula for $nb_n-(n-1)b_{n-1}$.

\begin{thm}
For $n\in\mathbb{N}$ and $0\le\ell\le n-1$, the combinatorial identity
\begin{equation}\label{G(mk)-Stirl}
\sum_{j=0}^{n}\frac{(-1)^j}{(n-j)!}G(j,\ell)=0
\end{equation}
is valid.
\par
For $k,n\in\mathbb{N}_0$, the partial Bell polynomials $\bell_{2n,k}$ can be computed by
\begin{multline}\label{Qi-Bell-Formula}
\bell_{2n,k}\biggl(0,-\frac{2!}{3}, 0, \frac{4!}{5}, \dotsc, \frac{1+(-1)^{2n-k+1}}{2}(-1)^{(2n-k+1)/2}\frac{(2n-k+1)!}{2n-k+2}\biggr)\\
=(-1)^{n-k}(2n)!\sum_{j=0}^{k} \frac{(-1)^j}{(k-j)!}G(j,n).
\end{multline}
\par
For $n\in\mathbb{N}$, the sequence $b_n$ can be explicitly computed by
\begin{equation}\label{b(n)-2sums}
b_n=(2n)!\sum_{j=1}^{2n}\frac{(-1)^{j+1}}{(2n-j)!j}G(j,n).
\end{equation}
\end{thm}

\begin{proof}
The proof proceeds in several steps.
\par
\paragraph{\indent\it Step~1: reduction to Bell polynomials.}
Let $u(x)=\frac{\arctan x}{x}$, so that $u(x)\to1$ as $x\to0$.
By~\eqref{b-n-gen-f-arctan}, it follows that
\[
(-1)^n(2n)!b_n=\lim_{x\to0}\frac{\operatorname{d}^{2n}}{\td x^{2n}}\ln u(x).
\] 
The Fa\`a di Bruno formula~\eqref{Bruno-Bell-Polynomial} with $f=\ln$ and $g=u$, together with $\ln^{(k)}(u)=\frac{(-1)^{k-1}(k-1)!}{u^k}$ and $\bell_{2n,0}=0$, gives
\begin{equation}\label{bn(2n)-1}
(-1)^n(2n)!b_n=\sum_{k=1}^{2n}(-1)^{k-1}(k-1)!\bell_{2n,k}(x_1,x_2,\dotsc,x_{2n-k+1}),
\end{equation}
where $x_m=\lim_{x\to0}u^{(m)}(x)$. From
\begin{equation*}
u(x)=\sum_{k=0}^{\infty}(-1)^k\frac{x^{2k}}{2k+1}, \quad |x|<1,
\end{equation*}
we have 
\[
u^{(m)}(x)=\sum_{k=0}^{\infty}(-1)^k\fall{2k}{m}\frac{x^{2k-m}}{2k+1}, 
\]
and therefore
\begin{equation*}
x_m=\begin{dcases}
0, & m=2\ell-1;\\
(-1)^\ell\frac{(2\ell)!}{2\ell+1}, & m=2\ell
\end{dcases}
\end{equation*}
for $\ell\in\mathbb{N}$.
Thus $(x_1,x_2,x_3,x_4,\dotsc)=\bigl(0,-\frac{2!}{3},0,\frac{4!}{5},\dotsc\bigr)$, which are the arguments appearing in~\eqref{Qi-Bell-Formula}.
\par
\paragraph{\indent\it Step~2: a formula for $\bell_{2n,k}$.}
Put 
\[
\varphi(z)=\sum_{m=0}^{\infty}\frac{x_{m+1}}{(m+1)!}z^m.
\] 
Then
\[
\sum_{m=1}^{\infty}x_m\frac{z^m}{m!}=z\varphi(z).
\]
Comparing the coefficients of $z^{N+k}$ in~\eqref{113-final-formula} gives 
\[
\bell_{N+k,k}(x_1,\dotsc,x_{N+1})=\binom{N+k}{k}\bigl(\varphi^k\bigr)^{(N)}(0). 
\]
Taking $N=2n-k$, we obtain
\begin{equation*}
\bell_{2n,k}(x_1,x_2,\dotsc,x_{2n-k+1})=\binom{2n}{k}\lim_{z\to0}\frac{\operatorname{d}^{2n-k}}{\td z^{2n-k}}\bigl[\varphi(z)\bigr]^k,\quad 0\le k\le 2n.
\end{equation*}
Since $x_{2j-1}=0$ and $x_{2j}=(-1)^j\frac{(2j)!}{2j+1}$, we have
\begin{equation*}
\varphi(z)=\sum_{j=1}^{\infty}(-1)^j\frac{z^{2j-1}}{2j+1}=\frac1z\biggl(\frac{\arctan z}{z}-1\biggr)=-\frac z3+\frac{z^3}5-\dotsm.
\end{equation*}
\par
\paragraph{\indent\it Step~3: expansion of $\varphi^k$ and the identity~\eqref{G(mk)-Stirl}.}
By the binomial theorem and the series expansion~\eqref{Grishaev-Sazonov-ITSF-1}, it follows that
\begin{multline*}
[\varphi(z)]^k=\frac{1}{z^k}\sum_{j=0}^{k}(-1)^{k-j}\binom{k}{j}\biggl(\frac{\arctan z}{z}\biggr)^j
\\
=\frac{(-1)^k}{z^k}\sum_{\ell=0}^{\infty}(-1)^\ell\Biggl[\sum_{j=0}^{k}(-1)^j j!\binom{k}{j} G(j,\ell)\Biggr]z^{2\ell}.
\end{multline*}
On the other hand, $\varphi(z)=-\frac z3+O(z^3)$, so $[\varphi(z)]^k=O(z^k)$, and the coefficients of $z^{2\ell-k}$ with $2\ell<2k$ must vanish. Using $j!\binom kj=\frac{k!}{(k-j)!}$, this is
\begin{equation*}
\sum_{j=0}^{k}\frac{(-1)^j}{(k-j)!}G(j,\ell)=0,\quad 0\le\ell\le k-1,\quad k\in\mathbb{N},
\end{equation*}
which is~\eqref{G(mk)-Stirl} (with $k$ renamed $n$). Consequently, with $\ell=k+i$, we obtain
\begin{equation}\label{power=k=eq}
[\varphi(z)]^k=\sum_{i=0}^{\infty}(-1)^{i}\Biggl[\sum_{j=0}^{k}(-1)^j j!\binom{k}{j} G(j,k+i)\Biggr]z^{k+2i},\quad k\in\mathbb{N}_0.
\end{equation}
\par
\paragraph{\indent\it Step~4: the formula~\eqref{Qi-Bell-Formula}.}
Let $0\le k\le n$. The term of~\eqref{power=k=eq} with $k+2i=2n-k$, that is, $i=n-k$, gives
\begin{equation*}
\lim_{z\to0}\frac{\operatorname{d}^{2n-k}}{\td z^{2n-k}}[\varphi(z)]^k=(-1)^{n-k}(2n-k)!\sum_{j=0}^{k}(-1)^j j!\binom{k}{j}G(j,n),
\end{equation*}
and hence, by Step~2 and $\binom{2n}{k}(2n-k)!j!\binom kj=\frac{(2n)!}{(k-j)!}$,
\begin{equation*}
\bell_{2n,k}(x_1,x_2,\dotsc,x_{2n-k+1})=(-1)^{n-k}(2n)!\sum_{j=0}^{k}\frac{(-1)^j}{(k-j)!}G(j,n).
\end{equation*}
This is~\eqref{Qi-Bell-Formula} for $0\le k\le n$. 
\par
For $k>n$, the left-hand side of~\eqref{Qi-Bell-Formula} vanishes, because $x_m=0$ for odd $m$ forces each part of a partition to be even, so a partition of $2n$ has at most $n$ parts; the right-hand side vanishes by~\eqref{G(mk)-Stirl} with $\ell=n\le k-1$.
\par
\paragraph{\indent\it Step~5: the formula~\eqref{b(n)-2sums}.}
Substituting the formula~\eqref{Qi-Bell-Formula} into~\eqref{bn(2n)-1}, and noting $(-1)^{k-1}(-1)^{n-k}=(-1)^{n-1}$, we obtain
\begin{equation*}
b_n=\sum_{k=1}^{2n}(k-1)!\sum_{j=0}^{k} \frac{(-1)^{j+1}}{(k-j)!}G(j,n),\quad n\in\mathbb{N}.
\end{equation*}
Since $G(0,n)=0$ for $n\in\mathbb{N}$, the term $j=0$ vanishes. Interchanging the order of summation gives
\begin{equation*}
b_n=\sum_{j=1}^{2n}(-1)^{j+1}G(j,n)\sum_{k=j}^{2n} \frac{(k-1)!}{(k-j)!}
\end{equation*}
and
\begin{equation*}
\sum_{k=j}^{2n} \frac{(k-1)!}{(k-j)!}
=\sum_{k=0}^{2n-j} \frac{(k+j-1)!}{k!}
=(j-1)!\sum_{k=0}^{2n-j} \binom{k+j-1}{j-1}
=\frac{1}{j}\frac{(2n)!}{(2n-j)!}
\end{equation*}
by the identity $\sum_{k=0}^{N}\binom{x+k}{x}=\binom{x+N+1}{x+1}$, see~\cite[pp.~27--28, (1.48) and~(1.49)]{Sprugnoli-Gould-2006}. This is~\eqref{b(n)-2sums}.
\end{proof}

\begin{rem}
In~\cite[Theorem~1.1]{AIMS-Math20210491.tex}, it was proved that the identity
\begin{equation*}
\bell_{2n+1,k}(0,x_2,0,x_4,\dotsc)=0, \quad n,k\in\mathbb{N}_0
\end{equation*}
is valid for any sequence $x_m\in\mathbb{C}$ such that $x_{2m-1}=0$ for $m\in\mathbb{N}$.
See also~\cite[Remark~3]{series-arccos-up-42.tex} and~\cite[p.~9347, Eq.~(2.7)]{integer2Bell2real.tex}.
In particular,
\begin{equation*}
\bell_{2n-1,k}\bigl(f'(0),f''(0),\dotsc,f^{(2n-k)}(0)\bigr)=0,\quad n\in\mathbb{N},\quad k\in\mathbb{N}_0,
\end{equation*}
whenever $f$ is an analytic even function on a symmetric interval about the origin.
\par
Similarly, if $x_{2m}=0$ for all $m\in\mathbb{N}$, then the left-hand side of~\eqref{113-final-formula} is an even function of $t$ for even $k$ and an odd function of $t$ for odd $k$, so $\bell_{N,k}=0$ whenever $N-k$ is odd. In particular,
\begin{equation*}
\bell_{2n,2m-1}(x_1,x_2,\dotsc,x_{2n-2m+1})=0, \quad n,m\in\mathbb{N}
\end{equation*}
and, if $f$ is an analytic odd function on a symmetric interval about the origin,
\[
\bell_{2n,2m-1}\bigl(f'(0),f''(0),\dotsc,f^{(2n-2m+2)}(0)\bigr)=0.
\]
\par
These identities can be applied to simplify or to understand more deeply some results on the partial Bell polynomials $\bell_{n,k}$ and their proofs in the works~\cite[Theorem~2; p.~454, Eqs.~(56) and~(57); p.~456, Eqs.~(58) and~(59)]{AADM-3164.tex}, \cite[Theorem~1.1]{AIMS-Math20210491.tex}, \cite[Theorems~4 and~5]{series-arccos-up-42.tex}, \cite[Lemma~1]{Note-Xu-Hoffman.tex}, \cite[Theorem~1.2]{Tan-Der-App-Thanks.tex}, \cite[Theorem~2.1]{Special-Bell2Euler.tex}, \cite[Lemma~2.5 and Remark~4.2]{KimKJK-JCAA-23-7.tex}, \cite{Char-Bell-Qi-MNN.tex}, \cite[Sections~1.6 and~1.11]{Bell-value-elem-funct.tex}, \cite[Theorem~7]{ser-sinc-real-pow-Qi-Taylor.tex}, \cite[Lemma~2.2 and Open Problem~4.1]{centr-bell-polyn.tex}, and~\cite[Lemmas~2.4 and~2.6]{Euler-No-3Sum.tex}.
\end{rem}

\begin{thm}
For $n\ge2$, we have
\begin{equation}\label{nb(n)-deff-eq}
nb_n-(n-1)b_{n-1}=\frac{1}{2}\sum_{j=1}^{2n} (-1)^jj! \binom{2n+1}{j+1}G(j,n)<0.
\end{equation}
\end{thm}

\begin{proof}
Taking $r=-1$ in~\eqref{Grishaev-Sazonov-ITSF-2}, where $\frac{(-r)_k}{k!}=1$, gives
\begin{equation}\label{recip=arctan-ser}
\frac{z}{\arctan z} =1+\sum_{n=1}^{\infty}(-1)^n\Biggl[\sum_{k=1}^{2n} \sum_{j=1}^{k}(-1)^j\binom{k}{j}j!G(j,n)\Biggr]z^{2n}, \quad |z|<1.
\end{equation}
Comparing~\eqref{recip=arctan-ser} with~\eqref{arctan-recipr} leads, for $n\ge2$, to
\begin{equation*}
2[nb_n-(n-1)b_{n-1}]=\sum_{k=1}^{2n} \sum_{j=1}^{k}(-1)^j\binom{k}{j}j!G(j,n)
=\sum_{j=1}^{2n} (-1)^jj!G(j,n)\sum_{k=j}^{2n}\binom{k}{j}
\end{equation*}
and $\sum_{k=j}^{2n}\binom{k}{j}=\binom{2n+1}{j+1}$ by~\cite[(1.48)]{Sprugnoli-Gould-2006}, that is, 
\[
\sum_{k=0}^{N}\binom{k+x}{r}=\binom{N+x+1}{r+1}-\binom{x}{r+1}
\]
with $r=j$, $x=0$, and $N=2n$.
This is the equality in~\eqref{nb(n)-deff-eq}; the inequality is Corollary~\ref{b(n)-positive-cor}.
\end{proof}

\section{Logarithmically complete monotonicity}\label{sec-logCM}

The positivity of $b_n$ has a natural interpretation in terms of monotonicity properties of functions. In this section, we first recall the relevant definitions and relations (Subsection~\ref{sec-lcm-background}). We then prove that $\frac{1}{\arctan x}$ is logarithmically completely monotonic on $(0,\infty)$ (Subsection~\ref{sec-lcm-arctan}) and that $\frac{\arctanh x}{x}$ is logarithmically absolutely monotonic on $(0,1)$ (Subsection~\ref{sec-labs-arctanh}).

\subsection{Definitions and background}\label{sec-lcm-background}

We recall the classes of functions involved and how they are related to each other.

Suppose that the function $h:I\to[0,\infty)$ is infinitely differentiable on an interval $I$. If
\begin{equation*}
0\le(-1)^{m-1}h^{(m-1)}(u)<\infty, \quad m\in\mathbb{N},
\end{equation*}
then we say that $h$ is a completely monotonic function on $I$. If $h:I\to(0,\infty)$ is positive and
\begin{equation*}
(-1)^m[\ln h(u)]^{(m)}\ge0, \quad m\in\mathbb{N},
\end{equation*}
then we say that $h$ is a logarithmically completely monotonic function on $I$. If $h:(0,\infty)\to[0,\infty)$ can be written in the form
\begin{equation*}
h(u)=\frac{a}u+b+\int_0^\infty\frac1{s+u}\td\mu(s),
\end{equation*}
where $a,b$ are non-negative constants and $\mu$ is a measure on $(0,\infty)$ such that $\int_0^\infty\frac1{1+s}\td\mu(s)<\infty$, then we say that $h$ is a Stieltjes function. We say that a function $h(x)$ is (logarithmically) absolutely monotonic on an interval $I$ if and only if $h(-x)$ is a (logarithmically) completely monotonic function on the interval $-I=\{-x:x\in I\}$.
For more details on definitions and properties of these kinds of functions, please refer to the papers~\cite{CBerg, absolute-mon-simp.tex, compmon2}, the chapters~\cite[Chapter~XIII]{mpf-1993} and~\cite[Chapter~IV]{widder}, and the monograph~\cite{Schilling-Song-Vondracek-2nd}.
Among these kinds of functions, there are the following relations:
\begin{enumerate}
\item
A function $h$ is completely monotonic on $(0,\infty)$ if and only if it is a Laplace transform, that is, there is a positive measure $\mu$ on $[0,\infty)$ such that $h(u)=\int_0^\infty \te^{-us}\td\mu(s)$.
\item
The set of all logarithmically completely monotonic functions is a strict subset of the set of all completely monotonic functions.
\item
The set of all Stieltjes functions is a strict subset of the set of all logarithmically completely monotonic functions on $(0,\infty)$.
\end{enumerate}
For details on these relations, please refer to~\cite{CBerg, absolute-mon-simp.tex, Qi-Agar-Surv-JIA.tex, compmon2, e-gam-rat-comp-mon, Schilling-Song-Vondracek-2nd} and closely related references therein.

\subsection{Logarithmically complete monotonicity of \texorpdfstring{$\frac1{\arctan x}$}{1/arctan x}}\label{sec-lcm-arctan}

In this subsection, we establish an integral representation of $\frac{1}{(1+x^2)\arctan x}$ as a Laplace transform, which yields the logarithmically complete monotonicity of $\frac{1}{\arctan x}$ and shows that it is not a Stieltjes function.

\begin{thm}\label{arctan-log-com-thm}
The integral representation
\begin{equation}\label{eq:gfin-s}
\frac{1}{(1+x^2)\arctan x} =2\int_0^\infty\biggl[\int_{0}^{\infty} \frac{1-\cos\bigl(s\coth\frac{v}{2}\bigr)}{\pi^2+v^2}\td v\biggr] \te^{-xs}\td s
\end{equation}
is valid for $x\in(0,\infty)$. Consequently, the function $\frac{1}{\arctan x}$ is logarithmically completely monotonic on $(0,\infty)$, but $\frac{1}{\arctan x}$ is not a Stieltjes function.
\end{thm}

\begin{proof}
In what follows, $\ln z$ denotes the principal value of the logarithm: $\ln z=\ln|z|+\ti \arg z$ with $\arg z\in(-\pi,\pi]$.
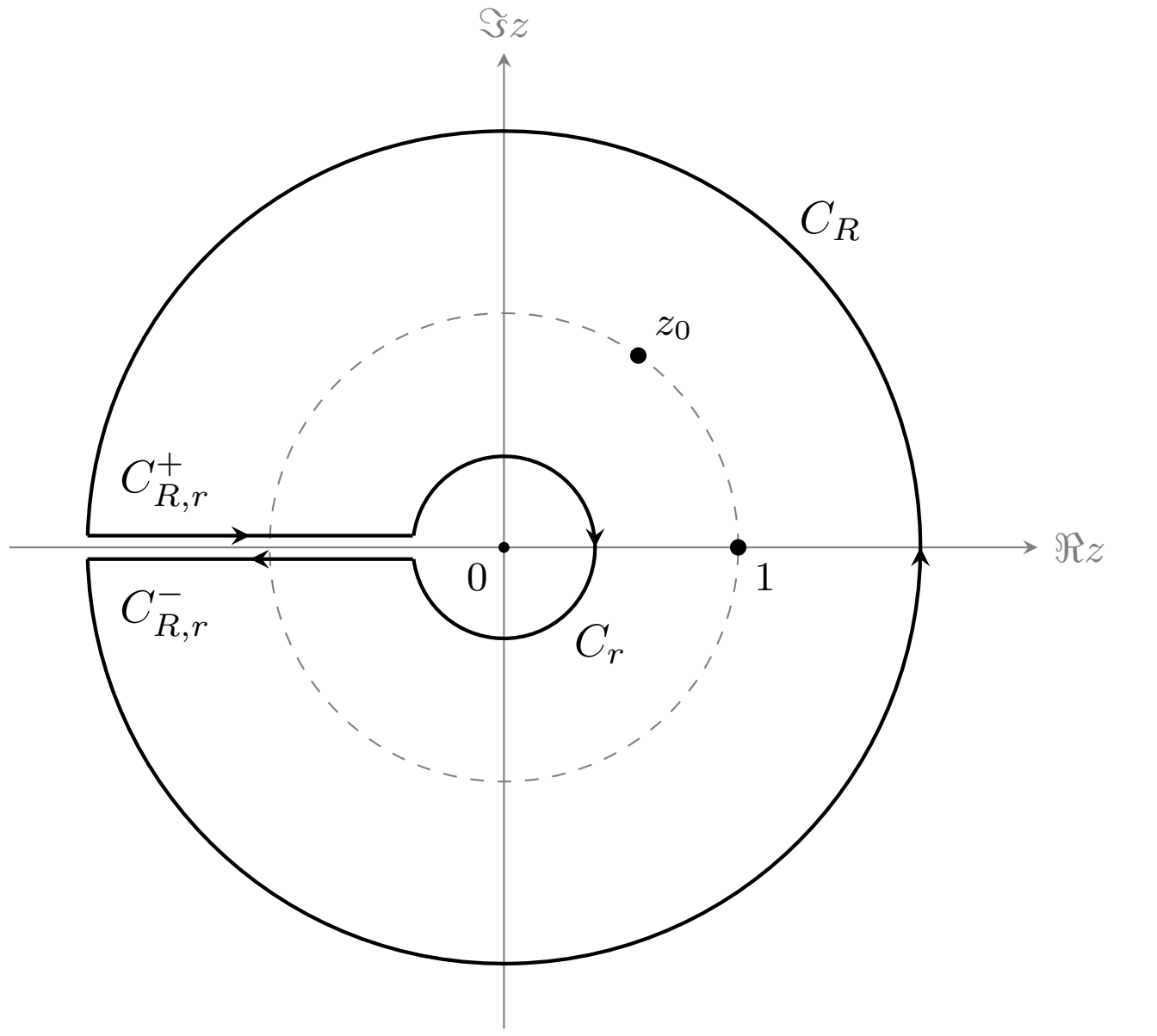
\begin{figure}[htbp]
\centering
\begin{tikzpicture}[>=stealth,
 mid/.style={decoration={markings,mark=at position 0.5 with {\arrow{>}}},postaction={decorate}}]
 \def\R{3.2}\def\r{0.7}\def\u{1.8}\def\e{0.09}
 \pgfmathsetmacro{\aR}{asin(\e/\R)}
 \pgfmathsetmacro{\ar}{asin(\e/\r)}
 \draw[->,gray] (-\R-0.6,0)--(\R+0.9,0) node[right]{\small$\Re z$};
 \draw[->,gray] (0,-\R-0.5)--(0,\R+0.6) node[above]{\small$\Im z$};
 \draw[dashed,gray] (0,0) circle (\u);
 \draw[thick,mid] ({\R*cos(-180+\aR)},{\R*sin(-180+\aR)}) arc (-180+\aR:180-\aR:\R);
 \draw[thick,mid] (-\R,\e)--(-\r,\e);
 \draw[thick,mid] ({\r*cos(180-\ar)},{\r*sin(180-\ar)}) arc (180-\ar:-180+\ar:\r);
 \draw[thick,mid] (-\r,-\e)--(-\R,-\e);
 \fill (\u,0) circle (1.8pt) node[below right]{\small$1$};
 \fill ({\u*cos(55)},{\u*sin(55)}) circle (1.8pt) node[above right]{\small$z_0$};
 \fill (0,0) circle (1.2pt) node[below left]{\small$0$};
 \node at (45:\R+0.35) {$C_R$};
 \node at (-45:\r+0.35) {$C_r$};
 \node at (-2.6,0.5) {$C^{+}_{R,r}$};
 \node at (-2.6,-0.5) {$C^{-}_{R,r}$};
\end{tikzpicture}
\caption{The contour $C=C_R\cup C_{R,r}^+\cup C_r\cup C_{R,r}^-$}
\label{fig:kontura}
\end{figure}
\par
Let us consider the integral $\int_{C} G(z)\td z$ along the positively oriented contour
\begin{equation*}
C=C_R\cup C_{R,r}^+\cup C_r\cup C_{R,r}^-
\end{equation*}
with $0<r<1<R$, as shown in Figure~\ref{fig:kontura}, where $C_R$ is the circle $|z|=R$ traversed counterclockwise, $C_{R,r}^{\pm}$ are the upper and lower edges of the cut along the negative real axis (the upper edge traversed from $-R$ to $-r$, the lower edge from $-r$ to $-R$), $C_r$ is the circle $|z|=r$ traversed clockwise, and
\begin{equation*}
G(z)=\frac{z+1}{z(z-z_0)\ln z}
\end{equation*}
with
\[
z_0=\frac{\ti-x_0}{\ti+x_0}=\frac{1-x_0^2+2\ti x_0}{1+x_0^2}
\] 
for $x_0>0$, so that $|z_0|=1$, $\Im z_0>0$, and $\ln z_0=2\ti\arctan x_0$.
\par
The points $1$ and $z_0$ are the only singularities of $G$ lying inside $C$.
By the residue theorem, we have
\begin{equation*}
\int_{C}G(z)\td z=2\pi\ti (\Res[G;z_0]+\Res[G;1]).
\end{equation*}
Since $z_0$ is a simple pole, using $1+z_0=\frac{2\ti}{\ti+x_0}$, we obtain
\begin{multline*}
\Res[G;z_0]=\frac{1+z_0}{z_0\ln z_0}
=\frac{\frac{2\ti}{\ti+x_0}}{\frac{\ti-x_0}{\ti+x_0}2\ti \arctan x_0}
\\
=\frac{1}{(\ti-x_0)\arctan x_0}
=-\frac{\ti+x_0}{(1+x_0^2)\arctan x_0}.
\end{multline*}
Similarly, since $\ln z\sim z-1$ as $z\to1$ and $1-z_0=\frac{2x_0}{\ti+x_0}$,
\begin{equation*}
\Res[G;1]=\lim_{z\to 1}\frac{(z-1)(1+z)}{z(z-z_0)\ln z}
=\frac{2}{1-z_0}=\frac{\ti+x_0}{x_0}.
\end{equation*}
As a result,
\begin{equation}\label{eq:g}
\frac{1}{(1+x_0^2)\arctan x_0}=\frac{1}{x_0}-\frac{1}{2\pi\ti(\ti+x_0)}\int_{C}G(z)\td z.
\end{equation}
\par
Now we calculate 
\[
\int_{C}G(z)\td z=\int_{C_R}G(z)\td z+\int_{C_r}G(z)\td z+\int_{C_{R,r}^+}G(z)\td z+\int_{C_{R,r}^-}G(z)\td z.
\]
On $|z|=R$ we have 
\[
|z+1|\le R+1,\quad |z-z_0|\ge R-1, \quad |\ln z|\ge\ln R,
\] 
while on $|z|=r$ we have 
\[
|z+1|\le1+r, \quad |z-z_0|\ge1-r,\quad |\ln z|\ge|\ln r|;
\] 
hence
\begin{equation*}
\biggl|\int_{C_R}G(z)\td z\biggr|\le 2\pi\frac{R+1}{(R-1)\ln R}
\quad\text{and}\quad
\biggl|\int_{C_r}G(z)\td z\biggr|\le2\pi\frac{1+r}{(1-r)|\ln r|},
\end{equation*}
which vanish as $R\to\infty$ and $r\to 0^+$, respectively.
\par
For $t<0$, the limits
\begin{equation*}
G^\pm(t)=\lim_{\substack{z\to t\\ \pm\Im z>0}}G(z)= \frac{t+1}{t(t-z_0)(\ln|t|\pm\pi\ti)}
\end{equation*}
exist, and the orientation of $C$ gives
\begin{equation*}
\int_{C_{R,r}^+}G(z)\td z+\int_{C_{R,r}^-}G(z)\td z
=\int_{-R}^{-r}\bigl[G^+(t)-G^-(t)\bigr]\td t.
\end{equation*}
Since 
\[
\frac{1}{\ln|t|+\pi\ti}-\frac{1}{\ln|t|-\pi\ti}=\frac{-2\pi\ti}{\ln^2|t|+\pi^2},
\] 
letting $R\to\infty$ and $r\to0^+$ and substituting $t\mapsto-t$ yield
\begin{align*}
\mathcal{I}&=\lim_{\substack{R\to\infty\\ r\to 0^+}}\int_{C}G(z)\td z
\\
&=\int_{-\infty}^0\frac{-2\pi\ti (t+1)}{t(t-z_0)\bigl(\ln^2|t|+\pi^2\bigr)}\td t\\
&= -2\pi\ti \int_{0}^\infty\frac{1-t}{t(t+z_0)\bigl(\pi^2+\ln^2t\bigr)}\td t.
\end{align*}
Using 
\[
t+z_0=\frac{x_0(t-1)+\ti(t+1)}{\ti+x_0}
\] 
and multiplying the numerator and denominator by the complex conjugate $x_0(t-1)-\ti(t+1)$, we arrive at
\begin{align*}
\mathcal{I}&=2\pi\ti (\ti+x_0)\int_{0}^\infty\frac{(1-t)^2x_0+\ti (1-t^2)}{t\bigl[x_0^2(1-t)^2+(1+t)^2\bigr]\bigl(\pi^2+\ln^2t\bigr)}\td t.
\end{align*}
Substituting this into~\eqref{eq:g} yields
\begin{equation*}
\frac{1}{(1+x_0^2)\arctan x_0}=\frac{1}{x_0}-\int_{0}^\infty\frac{(1-t)^2x_0+\ti (1-t^2)}{t[x_0^2(1-t)^2+(1+t)^2]\bigl(\pi^2+\ln^2t\bigr)}\td t.
\end{equation*}
Since the left-hand side is real, taking real parts gives
\begin{equation*}
\frac{1}{(1+x_0^2)\arctan x_0}=\frac{1}{x_0}-x_0\int_{0}^\infty\frac{(1-t)^2}{t[x_0^2(1-t)^2+(1+t)^2]\bigl(\pi^2+\ln^2t\bigr)}\td t.
\end{equation*}
\par
Since
\begin{equation*}
\int_0^{\infty}\frac{1}{t\bigl(\pi^2+\ln^2t\bigr)}\td t
=\int_{-\infty}^{\infty}\frac{1}{u^2+\pi^2}\td u
=1,
\end{equation*}
we can write 
\[
\frac1{x_0}=\frac1{x_0}\int_0^\infty\frac{x_0^2(1-t)^2+(1+t)^2}{t[x_0^2(1-t)^2+(1+t)^2]\bigl(\pi^2+\ln^2t\bigr)}\td t
\]
and obtain
\begin{equation*}
\frac{1}{(1+x_0^2)\arctan x_0}=\frac{1}{x_0}\int_{0}^\infty\frac{(1+t)^2}{t[x_0^2(1-t)^2+(1+t)^2]\bigl(\pi^2+\ln^2t\bigr)}\td t.
\end{equation*}
The integrand is invariant under the substitution $t\mapsto\frac1t$, so the integrals over $(0,1)$ and $(1,\infty)$ are equal, and
\begin{equation}\label{eq:g2}
\frac{1}{(1+x_0^2)\arctan x_0}=\frac{2}{x_0}\int_{0}^1\frac{(1+t)^2}{t[x_0^2(1-t)^2+(1+t)^2]\bigl(\pi^2+\ln^2t\bigr)}\td t.
\end{equation}
\par
For $a,b,x_0>0$, it is clear that
\begin{equation*}
\frac{1}{x_0(a^2x_0^2+b^2)}=\frac{1}{b^2}\biggl[\frac1{x_0}-\frac12\biggl(\frac{1}{x_0+b\ti/a}+\frac{1}{x_0-b\ti/a}\biggr)\biggr].
\end{equation*}
Using the formula
\begin{equation*}
\frac{\Gamma(z)}{s^z}=\int_{0}^{\infty} t^{z-1}\te^{-st}\td t,\quad \Re s>0,\quad \Re z>0,
\end{equation*}
see~\cite[p.~49]{Temme-96-book}, with $z=1$, where $\Gamma$ is the classical Euler gamma function, we derive
\begin{equation*}
\frac{1}{x_0\pm b\ti/a}=\int_0^\infty \te^{-(x_0\pm b\ti/a)s}\td s
\quad\text{and}\quad
\frac{1}{x_0}=\int_{0}^\infty \te^{-x_0s}\td s.
\end{equation*}
Consequently,
\begin{multline*}
\frac{1}{x_0(a^2x_0^2+b^2)}=\frac{1}{b^2}\int_{0}^\infty \biggl[1-\frac{\te^{-(b\ti/a)s}+\te^{(b\ti/a)s}}{2}\biggr]\te^{-x_0s}\td s
\\
=\frac{1}{b^2}\int_0^{\infty}\biggl(1-\cos\frac{bs}{a}\biggr)\te^{-x_0s}\td s.
\end{multline*}
Further setting $a=1-t$ and $b=1+t$ with $t\in(0,1)$ yields
\begin{equation*}
\frac{1}{x_0[x_0^2(1-t)^2+(1+t)^2]}=\frac{1}{(1+t)^2}\int_0^\infty \biggl[1-\cos\frac{(1+t)s}{1-t}\biggr]\te^{-x_0s}\td s.
\end{equation*}
Substituting this into~\eqref{eq:g2} and interchanging the order of integration (the integrand is non-negative) leads to
\begin{equation*}
\frac{1}{(1+x_0^2)\arctan x_0}
=2\int_0^\infty\biggl[\int_{0}^1\frac{1}{t\bigl(\pi^2+\ln^2t\bigr)} \biggl(1-\cos\frac{(1+t)s}{1-t}\biggr)\td t\biggr]\te^{-x_0s}\td s.
\end{equation*}
Since $x_0>0$ is arbitrary, this holds with $x_0$ replaced by $x$.
Using the substitution $t=\te^{-v}$, for which $\frac{1+t}{1-t}=\coth\frac v2$ and $\frac{\td t}{t}=-\td v$, we transform the inner integral into $\int_{0}^{\infty}\frac{1-\cos(s\coth\frac{v}{2})}{\pi^2+v^2}\td v$. This proves the integral representation~\eqref{eq:gfin-s}.
\par
From the Laplace transform~\eqref{eq:gfin-s}, the function $\frac{1}{(1+x^2)\arctan x}$ is completely monotonic on $(0,\infty)$, because the density is non-negative. Since
\begin{equation*}
-\biggl[\ln\frac{1}{\arctan x}\biggr]'=(\ln\arctan x)'=\frac{1}{(1+x^2)\arctan x},
\end{equation*}
the function $\frac{1}{\arctan x}$ is logarithmically completely monotonic on $(0,\infty)$.
\par
In~\cite[p.~457]{alzer-berg},~\cite[Remark~1.4]{CBerg},~\cite[Proposition~1.3]{Berg-Pedersen-JCAM-2001}, and~\cite[p.~507]{berg-pedersen-rocky}, it was recited that, if $f(x)\ne0$ is a Stieltjes function, so is the function $\frac1{xf(x)}$. If $\frac1{\arctan x}$ were a Stieltjes function, so would be $\frac{\arctan x}{x}$, and hence it would be completely monotonic on $(0,\infty)$. But
\begin{equation*}
\biggl(\frac{\arctan x}x\biggr)''=\frac{2 \bigl[\bigl(x^4+2x^2+1\bigr)\arctan x-x\bigl(2x^2+1\bigr)\bigr]}{x^3 (1+x^2)^2}
\end{equation*}
takes the negative value 
\[
\pi\sqrt3\,-\frac{45}8=-0.1836\dotsc
\] 
at the point $x=\frac{\sqrt3\,}3$. Hence $\frac1{\arctan x}$ is not a Stieltjes function.
\end{proof}

\begin{rem}
Theorem~\ref{arctan-log-com-thm} verifies a conjecture presented in~\cite[Eighth Problem]{Qi-Agar-Surv-JIA.tex} and later announced on the website \url{https://math.stackexchange.com/q/4247090} by the author. The conjecture itself was originally proposed by the author before 2010. Attempts to prove this conjecture were made in~\cite{Jovanovic-Treml-Kraguj-2025, Wrong-Proof}, but the correct result is provided solely in the arXiv preprint~\cite{Jovanovic-Treml-arXiv-2021}. The argument for Theorem~\ref{arctan-log-com-thm} presented here is a reformulation and modification of the proofs of~\cite[Theorems~2.1 and~2.2]{Jovanovic-Treml-arXiv-2021}.
\end{rem}

\subsection{Logarithmically absolute monotonicity of \texorpdfstring{$\frac{\arctanh x}{x}$}{arctanh(x)/x}}\label{sec-labs-arctanh}

Finally, we show that the positivity of $b_n$ established in Corollary~\ref{b(n)-positive-cor} implies the logarithmically absolute monotonicity of $\frac{\arctanh x}{x}$ on $(0,1)$.

\begin{thm}
The function $\frac{\arctanh x}{x}$ is logarithmically absolutely monotonic on $(0,1)$ and logarithmically completely monotonic on $(-1,0)$.
\end{thm}

\begin{proof}
By Theorem~\ref{arctanh-gen-bn-thm}, we have
\[
\ln\frac{\arctanh x}{x}=\sum_{n=1}^{\infty}b_nx^{2n}, \quad |x|<1,
\] 
and $b_n>0$ by Corollary~\ref{b(n)-positive-cor}. Differentiating termwise, for $m\in\mathbb{N}$ and $x\in(0,1)$,
\begin{equation*}
\biggl[\ln\frac{\arctanh x}{x}\biggr]^{(m)}=\sum_{2n\ge m}b_n\fall{2n}{m}x^{2n-m}>0,
\end{equation*}
so $\frac{\arctanh x}{x}$ is logarithmically absolutely monotonic on $(0,1)$. Since this function is even, $h(-x)=h(x)$, so by the definition it is logarithmically completely monotonic on $(-1,0)$.
\end{proof}

\section{Declarations}

\paragraph{\bf Author's Contributions}
The author contributed to all parts of the manuscript and read and approved the final manuscript.

\paragraph{\bf Funding}
The author was supported in part by the Youth Project of Hulunbuir City for Basic Research and Applied Basic Research (Grant No.~GH2024020) and by the Natural Science Foundation of Inner Mongolia Autonomous Region (Grant No.~2025QN01041).

\paragraph{\bf Informed Consent Statement}
Not applicable.

\paragraph{\bf Ethical Approval}
This research did not involve human participants or animals.

\paragraph{\bf Availability of Data and Material}
Data sharing is not applicable because no new data were created or analyzed in this study.

\paragraph{\bf Competing Interests}
The author declares no competing interests.

\paragraph{\bf Use of AI tools declaration}
The author declares that artificial intelligence (AI) tools were used in the preparation of this article.

\paragraph{\bf Institutional Review Board Statement}
Not applicable.

\paragraph{\bf Acknowledgements}
Not applicable.

\end{document}